\documentclass[oneside,a4paper]{article}
\usepackage{pp}
\usepackage{pictexwd,dcpic}
\usepackage[russian,english]{babel}
\usepackage[utf8]{inputenc}
\usepackage{perpage} 
\MakePerPage{footnote} 

\begin{document}

\title{Ramification conjecture and \\
Hirzebruch's property\\ of line arrangements}
\author{D. Panov and A. Petrunin}
\date{}
\maketitle

\nofootnote{A.~Petrunin was partially supported by NSF grant DMS 1309340.}
\nofootnote{D.~Panov is a Royal Society University Research Fellow.}
\begin{abstract}

The ramification of a polyhedral space is defined as the metric completion of the universal cover of its regular locus.

We consider mainly 
polyhedral spaces of two origins: quotients of Euclidean space by a discrete group of isometries and polyhedral metrics on $\CP^2$ with singularities at a collection of complex lines. 

In the former case we conjecture that quotient spaces always have a $\CAT[0]$ ramification and prove this in several cases.  In the latter case  we prove that the ramification is $\CAT[0]$ if the metric on $\CP^2$ is non-negatively curved. We deduce that complex line arrangements in $\CP^2$ studied by Hirzebruch have aspherical complement.
\end{abstract}

\section{Introduction}

The main objects of this article are Euclidean polyhedral spaces and their \emph {ramifications}. The ramification  of a polyhedral space is the metric completion of the universal cover of its regular locus. We are interested in the situation when the ramification is $\CAT[0]$. 

Two classes of polyhedral spaces that will play the most important role are quotients of $\RR^m$ by discrete isometric actions, and polyhedral K\"ahler manifolds; 
that is, polyhedral manifolds with a complex structure. 

\parbf{Quotients of $\RR^m$ and ramification conjecture.}
We start with the case of $\RR^m$ quotients where the ramification space admits an alternative description in terms of  arrangements of planes of (real) codimension $2$;
we will call such planes \emph{hyperlines}.

Consider a discrete isometric and orientation-preserving action $\Gamma\acts\RR^m$.
Denote by $\mathcal{L}_\Gamma$ the arrangement of all hyperlines which are fixed by at least one non-identical element in $\Gamma$. Define the ramification of $\Gamma\acts\RR^m$
(briefly $\Ram_\Gamma$)
as the universal cover of $\RR^m$  branching infinitely  along each hyperline in $\mathcal{L}_\Gamma$.

More precisely, if $\tilde W_\Gamma$ denotes the universal cover of
$$W_\Gamma=\RR^m\backslash
\left(\bigcup_{\ell\in\mathcal{L}_\Gamma}\ell\right)$$
equipped with the length metric induced from $\RR^m$
then $\Ram_\Gamma$ is the metric completion of $\tilde W_\Gamma$.

One of the main motivations of this paper is the following conjecture.

\begin{thm}{Ramification conjecture}\label{mainconjecture}
Let $\Gamma\acts\RR^m$ be a properly discontinuous isometric orientation-preserving action.
Then
\begin{enumerate}[a)]
\item $\Ram_\Gamma$ is a $\CAT[0]$ space.
\item The natural inclusion
$\tilde W_\Gamma\hookrightarrow \Ram_\Gamma$ is a homotopy equivalence.
\end{enumerate}

\end{thm}

Assume for an action $\Gamma\acts\RR^m$ that the ramification conjecture holds. 
Then since $\CAT[0]$ spaces are contactable,
$\tilde W_\Gamma$ is also contractible,
and so $W_\Gamma$ is aspherical.

The ramification conjecture generalizes a conjecture of Allcock  \cite[Conjecture 1.4]{allcock} on finite \emph{reflection groups} (recall that a reflection group is a discrete group generated by a set of reflections of a Euclidean space). 
Allcock considered the case of the action $\Gamma\acts\CC^m$ of a finite reflection group $\Gamma$ that complexifies the  orientation-reversing action of $\Gamma$ on $\RR^m$ generated by reflections. 
Allcock's conjecture is related to an earlier conjecture of Charney and Davis (see \cite[Conjecture 3]{charney-davis-95}) which in turn is motivated by a conjecture of Arnold, Pham and Thom on complex hyperplane  arrangements.

In the following theorem we collect the partial cases of the ramification conjecture which we can prove. 

\begin{thm}{Theorem}\label{thm:main}
The ramification conjecture holds in the following cases.
\begin{itemize}
\item[$(\mathrm{R}^+)$] If the action $\Gamma\acts\RR^m$ is the orientation-preserving index two subgroup of a reflection group.
\item[$(\mathbb{Z}_2)$] If $\Gamma$ is isomorphic to $\ZZ_2^k$.
\item[$(\RR^3)$]\label{3} If $m\le 3$.
\item[$(\CC^2)$]\label{4}
If $m=4$, and the action $\Gamma\acts\RR^4$
preserves a complex structure on $\RR^4$.
\end{itemize}

\end{thm}

The most involved case is $(\CC^2)$;
it is proved in Section~\ref{pkspaces} 
and relies on Theorem \ref{thm:3D-sphere},
which is the main technical result of this paper.

The proofs of other cases are simpler.
The case $(\mathrm{R}^+)$ follows from more general Proposition~\ref{prop:ramresh}.
In Section~\ref{3R}
we give two proofs of the case $(\RR^3)$: one is based on Theorem~\ref{3cone} and Zalgaller's theorem~\ref{thm:rigid-sphere-2D}
and the other on the case $(\mathrm{R}^+)$.

\begin{thm}{Corollary}\label{cor:braid}
Let $S_3\acts\CC^3$ be the action of symmetric group by permuting coordinates of $\CC^3$.
Then $\Ram_{S_3}$ is a $\CAT[0]$ space.
\end{thm}

The above corollary is deduced from the $(\CC^2)$-case of Theorem~\ref{thm:main}
since the action $S_3\acts\CC^3$ splits as a sum of an action on $\CC^2$ and a trivial action on $\CC^1$.
This corollary also follows from a result of Charney and Davis in \cite{charney-davis-93}.

\parbf{Polyhedral manifolds and Hirzebruch's question.}
Our study of ramifications of polyhedral manifolds sheds  some light on a question of Hirzebruch on complex line arrangements in $\CP^2$ asked in \cite{Hirzebruch}. To state this question recall the notion of complex reflection groups and arrangements. 

A {\it finite complex reflection} group is a group $\Gamma$ acting on $\CC^m$ by complex linear transformation generated by elements that fix a complex hyperplane in $\CC^m$. 
The arrangement of complex hyperplanes%
\footnote{that is, the set of hyperplanes fixed by at least one non-trivial element of $\Gamma$} 
 $\mathcal{L}_\Gamma$ in $\CC^m$ and its projectivization in $\CP^{m-1}$ are called {\it complex reflection arrangements}.

\begin{thm}{Hirzebruch's question, \cite{Hirzebruch}} Let  $\mathcal{L}$ be a complex line arrangement in $\CP^2$ consisting of $3\cdot n$ lines such that each line of $\mathcal{L}$ intersect others at exactly $n+1$ points. 
Is it true that $\mathcal{L}$ is a complex reflection arrangement?
\end{thm}
The above property will be called {\it Hirzebruch's property}.
Hirzebruch noticed that all complex reflection line arrangements in $\CP^2$ satisfy this property. These line arrangements consist of two infinite series and five exceptional examples. 
The infinite series are called $A_m^0$ ($m\ge 3$) and $A_m^3$ $(m\ge 2)$ and correspond to reflection groups $G(m,m,3)$ and $G(m,p,3)$ ($p<m$) from Shephard--Todd classification. 
Five  exceptional examples correspond to reflection groups $G_{23}, G_{24}, G_{25}, G_{26}, G_{27}$. 

Hirzebruch's question is still open, but we are able to prove the following.

\begin{thm}{Theorem}\label{cor:3n}
All line arrangements satisfying Hirzebruch's property have aspherical complements.
\end{thm}

Note that if the answer to Hirzebruch's question were positive, this theorem would follow from 
the work of Bessis \cite{bessis}. 
Bessis finished the proof of the old conjecture stating that complements to finite complex reflection arrangements are aspherical. 
Namely he proved this statement for the cases of groups  $G_{24}$, $G_{27}$, $G_{29}$, $G_{31}$, $G_{33}$ and $G_{34}$. 
As an immediate corollary of our theorem we get a new geometric proof of Bessis's theorem for the cases of groups  $G_{24}$ and $G_{27}$.

Theorem \ref{cor:3n} has a generalization to a larger class of arrangements, described in Corollary \ref{generalarrangement}. 
Note that on the one hand line arrangements with aspherical complements are quite rare, on the other hand no idea exists at the present of how to classify them.

\parbf{About the proof of Theorem \ref{cor:3n}.} 
It follows from \cite[Corollary 7.8]{panov}
that for any arrangement satisfying Hirzebruch's property except the union of three lines, there is a non-negatively curved polyhedral metric on $\CP^2$ with singularities at this arrangement. Hence to prove the theorem it is enough to show that the ramification of this polyhedral metric satisfies conditions a) and b) of Conjecture \ref{mainconjecture}.  Let us sketch how this is done.

Consider a $3$-dimensional pseudomanifold $\Sigma$ with a piecewise-spherical metric.
Define the \emph{singular locus} $\Sigma^{{\star}}$ of $\Sigma$
as the set of points in $\Sigma$ which do not admit a neighborhood isometric to an open domain in
the unit $3$-sphere.

Then the {\it ramification} of $\Sigma$
is defined as the  completion of the universal cover $\tilde {\Sigma}^\circ$ of the \emph{regular locus} $\Sigma^\circ=\Sigma\backslash\Sigma^{{\star}}$.
The obtained space will be denoted as $\Ram\Sigma$.

In Theorem \ref{thm:3D-sphere} we  characterize three-dimensional spherical polyhedral manifolds $\Sigma$
admitting an isometric $\RR^1$-action with geodesic orbits such that $\Ram\Sigma$ is $\CAT[1]$.
The key condition in Theorem \ref{thm:3D-sphere} is that all points in $\Sigma$  lie sufficiently close to the singular locus.

The existence of an $\RR^1$-action as above on $\Sigma$ is equivalent to the existence of
a complex structure on the Euclidean cone over $\Sigma$;
see Theorem~\ref{allPK}. 
The latter permits us to apply Theorem \ref{thm:3D-sphere} in the proof of Theorem \ref{cor:3n}
since any non-negatively curved polyhedral metric on $\CP^2$ has complex holonomy.
It follows then that the ramification of $\CP^2$ is locally $\CAT[0]$ and by an analogue of Cartan--Hadamard theorem it is globally $\CAT[0]$;
see Proposition~\ref{orbianalog}.

\section{More questions and observations}\label{sec:questions}

\parbf{Ramification of a polyhedral space.}
A Euclidean polyhedral space with
non-negative curvature in the sense of Alexandrov has to be a pseudomanifold, possibly with a non-empty boundary.

In fact, a stronger statement holds,
a Euclidean polyhedral space $\mathcal{P}$
has curvature bounded from below in the sense of Alexandrov
if and only if
its regular locus $\mathcal{P}^\circ$
is connected and convex in $\mathcal{P}$;
that is, any minimizing geodesic between points in $\mathcal{P}^\circ$ lies completely in $\mathcal{P}^\circ$ (compare \cite[Theorem 5]{milka}).

Recall that the {\it ramification} of $\mathcal{P}$
is defined as the  completion of the universal cover $\tilde {\mathcal{P}}^\circ$ of the \emph{regular locus} $\mathcal{P}^\circ$.
The Next question is intended to generalize the ramification conjecture 
to a wider setting that is not related to group actions.

\begin{thm}{Question}\label{quest:generalized}
Let $\mathcal{P}$ be a Euclidean polyhedral space. Suppose $\mathcal{P}$ has non-negative curvature in the sense of Alexandrov.
What additional conditions should be imposed on $\mathcal{P}$  to guarantee that $\Ram \mathcal{P}$ is a $\CAT[0]$ space
and the inclusion $\tilde{\mathcal{P}}^\circ\hookrightarrow\Ram \mathcal{P}$ is a homotopy equivalence?
\end{thm}

For a while we thought that no additional condition on $\mathcal{P}$ should be imposed;
that is, $\Ram \mathcal{P}$ is always a $\CAT[0]$ space.  
But then we found a counterexample in dimension 4 and higher;
see Theorem \ref{notCAT[1]}.

Nevertheless, Theorem~\ref{3cone} joined with Zalgaler's theorem \ref{thm:rigid-sphere-2D}
imply that no additional condition is needed if $\dim \mathcal{P}\le 3$.
Theorem~\ref{linearrangement} also gives an affirmative answer in a particular four-dimensional case.
The latter theorem is used to prove Theorem~\ref{cor:3n};
it also proves \cite[Conjecture 8.2]{panov}.

We do not know what conditions should be imposed in general if $\dim \mathcal{P}\ge 4$
but would like to formulate a conjecture in one interesting non-trivial case.

\begin{thm}{Conjecture} \label{planearrangement1} Let $\mathcal{P}$ be a Euclidean polyhedral space with non-negative curvature in the sense of Alexandrov.
Suppose  that $\mathcal{P}$ is homeomorphic to $\CP^m$ and its singularities form a complex hyperplane arrangement on $\CP^m$. 
Then $\Ram \mathcal{P}$ is  $\CAT[0]$ 
and the inclusion $\tilde{\mathcal{P}}^\circ\hookrightarrow\Ram \mathcal{P}$ is a homotopy equivalence.
\end{thm}

This conjecture holds for $m=2$ by Theorem \ref{linearrangement}. 
Existence of higher-dimensional examples 
of such polyhedral metrics on $\CP^m$ can be deduced from 
\cite{looijenga}.

\parbf{Two-convexity of the regular locus.}
The same argument as in \cite{panov-petrunin} shows that the regular locus $\mathcal{P}^\circ$ of a polyhedral space is two-convex,
that is, it satisfies the following property.

\emph{Assume $\Delta$ is a flat tetrahedron.
Then any locally isometric geodesic immersion in  $\mathcal{P}^\circ$
of three faces of $\Delta$ which agrees on three common edges  can be extended to a
locally isometric immersion $\Delta\looparrowright \mathcal{P}^\circ$.}

From the main result of Alexander, Berg and Bishop  in \cite{ABB},
it follows that every simply connected two-convex flat manifold with a smooth boundary is $\CAT[0]$. Therefore, if one could approximate  $(\Ram  \mathcal{P})^\circ$ by flat two-convex manifolds with smooth boundary, Alexander--Berg--Bishop theorem would imply that $\Ram  \mathcal{P}\z\in\CAT[0]$.

This looks as a nice plan to approach the problem, but it turns out that such a smoothing does not exist even for the action $\ZZ_2^2\acts\CC^2$ which changes the signs of the coordinates;
see the discussion after Proposition~5.3 in \cite{panov-petrunin}
or \emph{Two convexity} in \cite{petrunin-orthodox} for more details.

\parbf{Ramification around a subset.}
Given a subset $A$ in a metric space $X$,
define $\Ram_A X$ as  the completion of the universal cover of $X\backslash A$.
Then results of  Charney and Davis in \cite{charney-davis-93} imply the following.
\begin{enumerate}[(i)]
\item\label{charney-davis-S2} Let $x$, $y$ and $z$ be distinct points in $\mathbb{S}^2$.
Then
$$\Ram_{\{x,y,z\}}\mathbb{S}^2\in \CAT[1]$$
if and only if the triangle $[xyz]$ has perimeter $2\cdot\pi$.
In particular the points $x$, $y$ and $z$ lie on a great circle of $\mathbb{S}^2$.
\item\label{charney-davis-S3}
Let $X$, $Y$ and $Z$ be disjoint great circles in $\mathbb{S}^3$.
Then
$$\Ram_{X\cup Y\cup Z}\mathbb{S}^3\in\CAT[1]$$ if and only if $X$, $Y$ and $Z$ are fibers of the Hopf fibration $\mathbb{S}^3\to \mathbb{S}^2$ and their images $x,y,z\in \mathbb{S}^2$ satisfy condition (\ref{charney-davis-S2})\footnote{More precisely, the  quotient metric on the base $\mathbb{S}^2$ has curvature $4$, so $[xyz]$ should have perimeter $\pi$.}.
\end{enumerate}

The following two observations give a link between the above results and Question~\ref{quest:generalized}.

It turns out that if $\mathcal{P}_n$ is a sequence of two-dimensional spherical polyhedral spaces with exactly three singular points that approach $\mathbb{S}^2$ in the sense of Gromov--Hausdorff then
the limit position of singular points on $\mathbb{S}^2$ satisfies (\ref{charney-davis-S2}).

With a bit more work one can show a similar statement holds in the three-dimensional case.
More precisely, let $\mathcal{P}_n$ be a sequence of three-dimensional spherical polyhedral spaces with the singular locus formed by exactly three circles.
If $\mathcal{P}_n$ approaches $\mathbb{S}^3$ in the sense of Gromov--Hausdorff then the limit position of singular locus  satisfies (\ref{charney-davis-S3}).

We finish the discussion with one more conjecture.

\begin{thm}{Conjecture} Let $\mathcal H$ be
a complex hyperplane arrangement in $\CC^m$ . Then  $\Ram_{\mathcal H} \CC^m$ is $\CAT[0]$
if and only if the following condition holds.

Let $\ell$ be any complex hyperline\footnote{that is, an affine subspace of complex codimension 2}
that belongs to more than one complex hyperplane of $\mathcal H$.
Then for any complex hyperplane $h\subset \CC^m$ containing $\ell$ there is a hyperplane
$h'\in \mathcal H$ containing $\ell$ such that the angle between
$h'$ and $h$ is at most $\frac{\pi}{4}$.
\end{thm}

Note that all complex reflection hyperplane arrangements satisfy the conditions of this conjecture.
The two-dimensional version of this conjecture is Corollary \ref{linesinC2},
and the ``only if'' part follows from this corollary.
If this conjecture holds then, using the orbi-space version of Cartan--Hadamard theorem \ref{orbianalog}
and Allcock's lemma \ref{lem:allcock} in the same way as in the proof of Theorem \ref{thm:main},
one shows that the inclusion
$(\Ram_{\mathcal H}\CC^m)^\circ\hookrightarrow(\Ram_{\mathcal H}\CC^m)$ is a homotopy equivalence.
Hence this conjecture gives and alternative geometric approach to Bessis's result  \cite{bessis} on asphericity of complements to complex reflection arrangements.

\section{Preliminaries}\label{sec:prelim}

\parbf{Three types of ramifications.}
Recall that we consider three types of ramifications which are closely related:
for group actions, for polyhedral spaces and for subsets.

\begin{itemize}
\item Given a subset $A$ in a metric space $X$, we define $\Ram_A X$ as  the completion of the universal cover of $X\backslash A$. We assume here that $X\backslash A$ is connected.
\item Given a polyhedral space $\mathcal{P}$ (Euclidean, spherical or hyperbolic), the ramification $\Ram \mathcal{P}$ is defined as $\Ram_{A} \mathcal{P}$, where $A$ is the singular locus of $\mathcal{P}$.
\item Given an isometric and orientation-preserving action $\Gamma\acts\RR^m$, the ramification $\Ram_{\Gamma} \mathcal{P}$ can be defined as $\Ram (\RR^m/\Gamma)$;
this definition makes sense since $\RR^m/\Gamma$ is a polyhedral space.
\end{itemize}

\parbf{Curvature bounds for polyhedral spaces.}
\emph{A Euclidean polyhedral space} is a simplicial complex equipped with an intrinsic metric such that each simplex is isometric to a simplex in a Euclidean space.

\emph{A spherical polyhedral space} is a simplicial complex equipped with an intrinsic metric such that each simplex is isometric to a simplex in a unit sphere.

The link of any simplex in a polyhedral space
(Euclidean or spherical),
equipped with the angle metric forms a spherical polyhedral space.

The following two propositions give a more combinatorial description of polyhedral spaces with curvature bounded from below or above.

\begin{thm}{Proposition}\label{prop:poly-cbb}
An $m$-dimensional Euclidean (spherical) polyhedral space $\mathcal{P}$
has curvature $\ge 0$ (correspondingly $\ge 1$)
in the sense of Alexandrov
if and only if each of the following conditions holds.
\begin{enumerate}
\item The link of any $(m-1)$-simplex is isometric to the one-point space $\mathfrak{p}$ or $\mathbb{S}^0$; 
that is, the two-point space with distance $\pi$ between the distinct points.
\item The link of any $(m-2)$-simplex is isometric to a closed segment of length $\le \pi$ or a circle with length $\le2\cdot\pi$.
\item The link of any $k$-simplex with $k\le m-2$ is connected.
\end{enumerate}

\end{thm}

\begin{thm}{Corollary}
The simplicial complex of any polyhedral space $\mathcal{P}$ with a lower curvature bound is a pseudomanifold.
\end{thm}

The following proposition
follows from Cartan--Hadamard theorem
and its analogue is proved by Bowditch in \cite{bowditch};
see also \cite{akp}, where both theorems are proved nicely.

\begin{thm}{Proposition}\label{prop:poly-cba}
A  polyhedral space $\mathcal{P}$ is a $\CAT[0]$ space
if and only if $\mathcal{P}$ is simply connected
and the link of each vertex is a $\CAT[1]$ space.

A spherical polyhedral space $\mathcal{P}$
is a $\CAT[1]$ space
if and only if  the link of each vertex of $\mathcal{P}$ is a $\CAT[1]$ space
and any closed curve of length $<2\cdot\pi$ in $\mathcal{P}$ is null-homotopic in the class of curves of length $<2\cdot\pi$.
\end{thm}

We say that a polyhedral space has \emph{finite shapes} if the number of isometry types of simplices that compose it is finite.
The following proposition is proved in  \cite[II. 4.17]{bridson-haefliger}.

\begin{thm}{Proposition}\label{geodesic2pi} Let $\mathcal P$ be a Euclidean  (spherical) polyhedral space with finite shapes and suppose that $\mathcal P$ has curvature $\le 0$ (or $\le 1$ correspondingly).
If $\mathcal P$ is not a $\CAT[0]$ (correspondingly, not $\CAT[1]$) space, 
then $\mathcal P$ contains an isometrically embedded
circle (correspondingly, a circle with length smaller than $2{\cdot}\pi$).
\end{thm}

\parbf{Spherical polyhedral metrics on $\bm{\mathbb{S}^2}$.}
The following theorem appears as an intermediate statement
in Zalgaller's proof of rigidity of spherical polygons;
see \cite{zalgaller}.

\begin{thm}{Zalgaller's theorem}\label{thm:rigid-sphere-2D}
Let $\Sigma$ be a spherical polyhedral space homeomorphic to the $2$-sphere
and with curvature $\ge 1$ in the sense of Alexandrov.
Assume that there is a point $z\in \Sigma$ such that all singular points lie at the distance $>\tfrac\pi2$ from $z$.
Then $\Sigma$ is isometric to the standard sphere.

\end{thm}

\parit{A sketch of Zalgaller's proof.}
We apply an induction on the number $n$ of singular points.
The base case $n=1$ is trivial.
To do the induction step
choose two singular points $p,q\in\Sigma$,
cut $\Sigma$ along a geodesic $[pq]$
and patch the hole so that the obtained new polyhedron $\Sigma'$ has curvature $\ge 1$.
The patch is obtained by doubling%
\footnote{Given a metric length space $X$ with a closed subset $A\subset X$,
the \emph{doubling} of $X$ across $A$ is obtained by gluing two copies of $X$ along $A$.}
 a convex spherical triangle across two sides.
For a unique choice of triangle, the points $p$ and $q$ become regular in $\Sigma'$
and exactly one new singular point appears in the patch.%
\footnote{This patch construction was introduced by Alexandrov,
the earliest reference we found is
\cite[VI, \S7]{alexandrov1948}.}
In this way, the case with $n$ singular points is reduced to the case with $n-1$ singular points.
\qeds

\parbf{A test for homotopy equivalence.}
The following lemma is a slight modification of Lemma 6.2 in \cite{allcock};
the proofs of these lemmas are almost identical.

\begin{thm}{Allcock's lemma}\label{lem:allcock}
Let $\mathcal{S}$ be an $m$-dimensional \emph{pure}%
\footnote{that is, each simplex in $\mathcal{S}$ forms a face in an $m$-dimensional simplex.} simplicial complex.

Let $K$ be a subcomplex in $\mathcal{S}$ of codimension $\ge 1$;
set $W=\mathcal{S}\backslash K$.
Assume that the link in  $\mathcal{S}$  of any simplex in $K$ is contractible.
Then the inclusion map $W\hookrightarrow\mathcal{S}$ is a homotopy equivalence.
\end{thm}

\parit{Proof.}
Denote by $K_n$ the $n$-skeleton of $K$;
set $W_n=\mathcal{S}\backslash K_n$
and
set $K_{-1}=\emptyset$.

For each $n\in \{0,1,\dots,m-1\}$
we will construct a homotopy
$$F_n\:[0,1]\times W_{n-1}\to W_{n-1}$$
of the identity map $\id_{W_{n-1}}$
into a map with the target in $W_{n}$.

Note that $W_{m-1}=W$ and $W_{-1}=\mathcal{S}$.
Therefore joining all the homotopies $F_n$,
we construct a  homotopy of the identity map on ${\mathcal{S}}$ into a map with the target in $W$.
Therefore the lemma follows once we construct  $F_n$ for all $n$.

\parit{Existence of $F_n$.}
Note that each open $n$-dimensional simplex $\Delta$ in $\mathcal{S}$ admits a closed neighbourhood $N_\Delta$ in $W_{n-1}$ which is homeomorphic to
$$\Delta\times (\mathfrak{p}\star\Link\Delta),$$
where $\Link\Delta$ denotes the link of $\Delta$,
$\mathfrak{p}$ denotes a one-point complex, and ${\star}$ denotes the join.
Moreover, we can assume that $\Delta$ lies in $N_\Delta=\Delta\times (\mathfrak{p}*\Link\Delta)$ as $\Delta\times \mathfrak{p}$
and $N_\Delta\cap N_{\Delta'}=\emptyset$ for any two open $(n-1)$-dimensional simplexes $\Delta$ and $\Delta'$ in $\mathcal{S}$.

Note that if $\Link\Delta$ is contractible
then  $\Link\Delta$ is a strict deformation retract of $\mathfrak{p}\star\Link\Delta$.
It follows that for any $(n-1)$-dimensional simplex $\Delta$ in $K$,
the relative boundary  $\partial_{W_{n-1}}N_\Delta$
is a deformation retract of $N_\Delta$.
Clearly $\partial_{W_{n-1}}N_\Delta\subset W_n$.
Hence the existence of $F_n$ follows.
\qeds

\parbf{An orbi-space version
of Cartan--Hadamard theorem.}

\begin{thm}{Proposition}\label{orbianalog}
Let $\mathcal P$ be a polyhedral pseudomanifold.
Suppose that for any point $x\in \mathcal P$ the ramification of the cone at $x$ is $\CAT[0]$.
Then
\begin{enumerate}[(i)]
\item\label{orbianalog:i} $\Ram \mathcal P$ is $\CAT[0]$;
\item\label{orbianalog:ii} for any $y\in \Ram \mathcal P$ that projects to $x\in \mathcal P$
the cone at $y$ is isometric to the ramification of the cone at $x$.
\end{enumerate}
\end{thm}

The proposition can be proved along the same lines as Cartan--Hadamard theorem;
see for example \cite{akp}.

A closely related statement was rigorously proved by Haefliger in \cite{haefliger};
he showed that if the charts of an orbi-space are $\CAT[0]$
then its universal orbi-cover is $\CAT[0]$.
Haefliger's definition of orbi-space restricts only to finite isotropy groups, but the above proposition requires only minor modifications of  Haefliger's  proof.

\parbf{Polyhedral K\"ahler manifolds.}
Let us recall some definitions and results
from \cite{panov}.
We will restrict our consideration to the case of non-negatively curved polyhedra.

\begin{thm}{Definition}\label{PKdefinition} Let $\mathcal{P}$ be an orientable non-negatively curved Euclidean
polyhedral manifold of dimension $2{\cdot}n$.
We say that  $\mathcal{P}$ is
\emph{polyhedral K\"ahler} if the holonomy of the metric
on $\mathcal{P}^\circ$ belongs to
$\mathrm{U}(n)< \SO(2{\cdot}n)$.

In the case  when $\mathcal{P}$ is a metric cone
piecewise linearly isomorphic to $\mathbb R^{2\cdot n}$ we call it a \emph{polyhedral K\"ahler cone}.
\end{thm}

Recall that from a result of Cheeger (see \cite{cheeger} and \cite[Proposition 2.3]{panov})
it follows that the metric of an orientable simply connected non-negatively curved polyhedral
compact $4$-manifold not homeomorphic to $\mathbb{S}^4$ has unitary holonomy.
Moreover in the case when the (unitary) holonomy is irreducible,
the manifold has to be diffeomorphic to $\CP^2$. Metric singularities
form a collection of complex curves on $\CP^2$;
see \cite{panov} for the details.

The following theorem summarizes some results on non-negatively
curved four-dimensional polyhedral K\"ahler cones proven in \cite[Theorems 1.5, 1.7 and 1.8]{panov}.

\begin{thm}{Theorem}\label{allPK}
Let $\mathcal C^4$ be a non-negatively curved polyhedral K\"ahler cone and let $\Sigma$
be the unit sphere of this cone centered at its tip.

\begin{enumerate}[(a)]
\item There is a canonical isometric $\RR$-action on $\mathcal C$
such that its orbits on $\Sigma$ are geodesics. This action is generated
by the vector field $J(r\frac{\partial}{\partial r})$ in the non-singular part
of $\mathcal C$, where $J$ is the complex structure on $\mathcal C$ and $r\frac{\partial}{\partial r}$
is the radial vector field on $\mathcal C$.
\item If the metric singularities of the cone are topologically equivalent to
a collection of $n\ge 3$ complex lines in $\CC^2$, then the action
of $\RR$ on $\Sigma$ factors through $\mathbb{S}^1$ and the map
$\Sigma\to \Sigma/\mathbb{S}^1$ is the Hopf fibration.
\item\label{allPK:3} If the metric singularities of the cone are topologically equivalent to
a union of two complex lines, then the cone splits as a metric product
of two two-dimensional cones.
\end{enumerate}
\end{thm}

\parbf{Reshetnyak gluing theorem.}
Let us recall the formulation of Reshetnyak gluing theorem which will be used in the proof of  Proposition \ref{prop:ramresh}.

\begin{thm}{Theorem}\label{thm:gluing}
Suppose that
$\spc{U}_1, \spc{U}_2$ are $\CAT[\kappa]$ spaces\footnote{We always assume that $\CAT[\kappa]$ spaces are complete.} with
closed convex subsets $A_i\subset\spc{U}_i$
which admit an isometry $\iota\:A_1\to A_2$.
Let us define a new space $W$ by gluing 
$\spc{U}_1$ and  $\spc{U}_2$ along the isometry $\iota$;
that is,  consider the new space
\[\spc{W}=\spc{U}_1\sqcup_\sim\spc{U}_2\]
where the equivalence relation $\sim$ is defined by $a\sim \iota(a)$ with the induced length metric. Then the following holds.

The space $\spc{W}$ is $\CAT[\kappa]$.
Moreover, both canonical mappings
$\tau_i\:\spc{U}_i\to\spc{W}$ are distance preserving,
and the images $\tau_i(\spc{U}_i)$ are convex subsets in $\spc{W}$.
\end{thm}

The following corollary is proved by repeated application of Reshetnyak's theorem.

\begin{thm}{Corollary}\label{cor:resh}
Let $S$ be a finite tree.
Assume a convex Euclidean (or spherical) polyhedron $Q_\nu$
corresponds to each node $\nu$ in $S$
and for each edge $[\nu\mu]$ in $S$
there is an isometry $\iota_{\mu\nu}$
from a facet\footnote{A facet is a face of codimension 1} $F\subset Q_\nu$ to a facet $F'\subset Q_\mu$.

Then the space obtained by gluing all the polyhedra $Q_\nu$ along the isometries $\iota_{\mu\nu}$ forms a $\CAT[0]$ space (correspondingly a $\CAT[1]$ space).
\end{thm}

\parbf{Flag complexes.}

\begin{thm}{Definition}
A simplicial complex $\mathcal{S}$ is \emph{flag} if whenever $\{v_0,\z\dots,v_k\}$
is a set of distinct vertices which are pairwise joined by edges, then $\{v_0,\dots,v_k\}$
spans a $k$-simplex in $\mathcal{S}$.
\end{thm}

Note that every flag complex is determined by its 1-skeleton.

Spherical polyhedral $\CAT[1]$ spaces glued from right-angled simplices
admit the following combinatorial characterization discovered by Gromov \cite[p. 122]{gromov-hyp-group}.

\begin{thm}{Theorem}
 A piecewise-spherical simplicial complex made of right-angled simplices is a $\CAT[1]$ space if and only if it is a flag complex.
\end{thm}

\section{On the reflection groups}

If the singular locus of a polyhedral space $\mathcal{P}$ coincides with its $(m-2)$-skeleton
then $\mathcal{P}^\circ$ has the homotopy type of a graph 
(its vertices correspond to the centers of $m$-simplices of $\mathcal{P}$).
We will show that in this case the ramification conjecture can be proven by applying Reshetnyak gluing theorem recursively.

We will prove the following stronger statement.

\begin{thm}{Proposition}\label{prop:ramresh}
Assume $\mathcal{P}$ is an $m$-dimensional polyhedral space which
admits a subdivision into closed sets $\{Q_i\}$
such that each $Q_i$ with the induced metric is isometric to a convex $m$-dimensional polyhedron and each face  of dimension $m-2$ of each polyhedron $Q_i$ belongs to the singular locus of $\mathcal{P}$.
Then $\Ram \mathcal{P}\in\CAT[0]$.
\end{thm}

Note that Theorem~\ref{thm:main}$(\mathrm{R}^+)$
follows directly from the above proposition.
Also the condition in the above proposition holds if $\mathcal{P}$ is isometric to the boundary of a convex polyhedron in Euclidean space
and in particular, by Alexandrov's theorem it holds if $\mathcal{P}$ is homeomorphic to $\mathbb{S}^2$.

\parit{Proof of Proposition \ref{prop:ramresh}.}
In the subdivision of $\mathcal{P}$ into $Q_i$,
color all the facets in different colors.
Consider the graph $\Gamma$ with a node for each $Q_i$, where two nodes are connected by an edge if the correspondent polyhedra have a common facet. Color each edge of $\Gamma$ in the color of the corresponding facet.

Denote by $\tilde\Gamma$ the universal cover of $\Gamma$.
Note that $\tilde\Gamma$ has to be a tree.

For each node $\nu$ of $\tilde\Gamma$,  prepare a copy of $Q_i$ which corresponds to the projection of $\nu$ in $\Gamma$.

Note that  the space $\Ram\mathcal{P}$
can be obtained by gluing the prepared copies.
Two copies should be glued along two facets of the same color $z$
if the nodes corresponding to these copies are connected in $\tilde\Gamma$
by an edge of color $z$.

Given a finite subtree $S$ of $\tilde \Gamma$
consider the subset $Q_S\subset \Ram\mathcal{P}$
formed by all the copies of $Q_i$ corresponding to the nodes of $S$.

Note that $Q_S$ is a convex subset of $\Ram\mathcal{P}$.
Indeed, if a path between points of $Q_S$ escapes from $Q_S$,
it has to cross the boundary $\partial Q_S$
at the same facet twice, say at the points $x$ and $y$ in 
a facet $F\subset\partial Q_S$.
Further note that the natural projection $\Ram\mathcal{P}\to \mathcal{P}$ is a short map which is distance preserving on $F$.
Therefore there is a unique geodesic from $x$ to $y$ and it lies in $F$.
In particular, geodesic with ends in $Q_S$ cannot escape from $Q_S$;
in other words $Q_S$ is convex.

Finally, by Corollary~\ref{cor:resh},
the subspace $Q_S$ is  $\CAT[0]$ for any finite subtree $S$.
Clearly, for every triangle $\triangle$ in $\Ram\mathcal{P}$
there is a finite subtree $S$ such that $Q_S\supset\triangle$.
Therefore the $\CAT[0]$ comparison holds for any geodesic triangle in $\Ram\mathcal{P}$.
\qeds

\section{Case $(\mathbb{Z}_2)$}

In this section we reduce the case $(\mathbb{Z}_2)$ of Theorem \ref{thm:main} to  the case $(\mathrm{R}^+)$.

\parit{Proof of Theorem \ref{thm:main}; case $(\mathbb{Z}_2)$.}
Every orientation-preserving action of a group $\ZZ_2^k$ on $\RR^m$ arises as
the action of a subgroup of the group $\ZZ_2^m$ generated by reflections in coordinate  hyperplanes.
By the definition of ramification, 
we can assume that the action of $\ZZ_2^k$ 
is generated by reflections in hyperlines.

Let us write $i\sim j$ if $i=j$ or  the reflection in
the hyperline $x_i=x_j=0$ belongs to $\Gamma$.
Note that $\sim$ is an equivalence relation.

It follows that $\RR^m/\Gamma$ splits as a direct product of the subspaces corresponding to the coordinate subspaces of $\RR^m$ for each equivalence relation.

Finally, for each of the factors in this splitting,
the statement holds by Theorem~\ref{thm:main} $(\mathrm{R}^+)$.
\qeds

\section{2-spaces}\label{ramsphere}

\begin{thm}{Definition} An $m$-dimensional spherical polyhedral space  $\Sigma$ is called \emph{$\alpha$-extendable}
if for any $\eps>0$,
every isometric immersion into $\Sigma$
of a ball of radius $\alpha+\epsilon$
from $\mathbb{S}^m$ extends to an isometric immersion of the whole $\mathbb{S}^m$.

In other words
$\Sigma$ is $\alpha$-extendable if
either
the distance from any point $x\in \Sigma$
to its singular locus $\Sigma^{{\star}}$ is at most $\alpha$
or $\Sigma$ is a space form.
\end{thm}

\begin{thm}{Theorem}\label{thm:2D-sphere}
Let $\Sigma$ be a two-dimensional spherical polyhedral manifold.
 Then $\Ram\Sigma$ is $\CAT[1]$ if and only if $\Sigma$ is $\frac{\pi}{2}$-extendable.
\end{thm}

\parit{Proof.}
Note that if $\Sigma^{{\star}}=\emptyset$ then
$\Sigma$ is a spherical space form.
So we assume that $\Sigma^{{\star}}\ne\emptyset$.

Let us show that in the case when $\Sigma$ is $\frac{\pi}{2}$-extendable one can
decompose $\Sigma$ into a collection of convex spherical
polygons with vertices in $\Sigma^{{\star}}$.
The proof is almost identical to the
proof of \cite[Proposition 3.1]{thurston},
so we just recall the construction.

In the case where $\Sigma^{{\star}}$ consists of two points,
$\Sigma$ can be decomposed into a collection of two-gons.
It remains to consider the case when $\Sigma^{{\star}}$ has at least three distinct points.

Consider the Voronoi decomposition of $\Sigma$ with respect to the points in $\Sigma^{{\star}}$.
The vertices of this decomposition
consist of points $x$ that have the following property.
If $D$ is the maximal open ball in $\Sigma^\circ=\Sigma\backslash \Sigma^{{\star}}$
with the center at $x$, then the radius of $D$ is at most $\tfrac\pi2$ and the convex hull of points in  $\partial D\cap \Sigma^{{\star}}$
contains $x$.
Note that such a convex hull is a convex spherical polygon $P(x)$ and $\Sigma$ is decomposed
into the union of $P(x)$ for various vertices $x$.

Consider finally the Euclidean cone over $\Sigma$ with the induced decomposition
into cones over spherical polygons. Applying Proposition \ref{prop:ramresh} to the cone we
see that its ramification is $\CAT[0]$. So $\Ram\Sigma\in \CAT[1]$ by Proposition \ref{prop:poly-cba}.
\qeds

The next result follows directly from Theorem \ref{thm:2D-sphere} and Proposition \ref{prop:poly-cba}.

\begin{thm}{Theorem}\label{3cone}
Let $\spc{Y}$ be a three-dimensional polyhedral cone.
Then $\Ram \spc{Y}$ is $\CAT[0]$
if and only if
$\spc{Y}$ satisfies one of the following conditions.
\begin{enumerate}
\item The singular locus $\mathcal{Y}^\star$ is formed by the tip or it is empty.
\item For any direction $v\in\spc{Y}$ there is a direction $w\in\mathcal{Y}^\star$ such that $\measuredangle(v,w)\le \tfrac\pi2$.
\end{enumerate}
\end{thm}
Indeed, the link $\Sigma$ of $\spc{Y}$ is a space form 
if and only if
$\mathcal{Y}^\star$ is formed by the tip or it is empty.
If $\mathcal{Y}^\star$ contains more than one point
the condition $2$ of this theorem means literally that $\Sigma$ is $\frac{\pi}{2}$-extendable.

\section{Case $(\RR^3)$}\label{3R}

Here we present two proofs of Theorem \ref{thm:main} case $(\RR^3)$: 
the first one is based on Theorem~\ref{3cone} and the second 
on Theorem \ref{thm:main} case $(\mathrm{R}^+)$.

In both of these proofs we assume that $\Gamma$ is finite.
The case when $\Gamma$ is infinite can be done the same way as the $\CC^2$ case;
see Section \ref{pkspaces}.

\parit{Proof 1.}
Since $\Gamma$ is finite, 
without loss of generality we may assume that $\Gamma$ fixes the origin.

By Zalgaller's theorem \ref{thm:rigid-sphere-2D}, the link of the origin in the 
quotient $\RR^3/\Gamma$ is $\tfrac\pi2$-extendable.
Applying Theorem \ref{3cone}, we get the result.
\qeds

\parit{Proof 2.}
By Theorem \ref{thm:main}
it is sufficient to prove that
$\Gamma$ is an index-two subgroup in a group $\Gamma_1$ generated by reflections in planes.

If $\Gamma$ fixes a line in $\RR^3$ then it is a cyclic group and it is an index-two subgroup of
a dihedral group.

Otherwise $\mathbb{S}^2/\Gamma$ is an orbifold with three orbi-points glued
from two copies of a Coxeter spherical triangle.
Such an orbifold has an involution
$\sigma$ such that $(\mathbb{S}^2/\Gamma)/\sigma$ is a Coxeter triangle $\Delta$. So $\Gamma_1$
is the group generated by reflections in the sides of $\Delta$.
\qeds

\section{3-spaces with a geodesic actions}

The following theorem is the main technical result.

\begin{thm}{Theorem}\label{thm:3D-sphere}
Let $\Sigma$ be a three-dimensional spherical polyhedral manifold.
Assume $\Sigma$ admits an isometric action of $\RR$ with geodesic orbits.

Then $\Ram\Sigma$ is $\CAT[1]$ if and only if  $\Sigma$ is $\frac\pi4$-extendable
or $\Ram\Sigma$ is the completion of the universal cover of $\mathbb{S}^3\setminus \mathbb{S}^1$.
\end{thm}

{\bf Example.} We will further apply this theorem to unit spheres of polyhedral
cones that are quotients of $\CC^2$ by a finite group of unitary isometries. The action
of $\RR$ in this case comes from the action on $\CC^2$ by multiplication
by complex units.

\medskip

The proof of
Theorem \ref{thm:3D-sphere} relies on several lemmas.
The following lemma
is spherical analogue of the theorem proved by German Pestov and Vla\-di\-mir Ionin in \cite{pestov-ionin};
a different proof via curve-shortening flow was given by Konstantin Pankrashkin in  \cite{pankrashkin};
see also \emph{The moon in the puddle} in
\cite{petrunin-orthodox}.

\begin{thm}{Drop lemma}\label{thickdrop}
Let $D$ be a disk with a metric of curvature $1$,
whose boundary consists of several smooth arcs  of curvature at most $\kappa$ that meet
at angles larger than $\pi$ at all points except at most one. 
Then:

\begin{enumerate}[(a)]
\item\label{thickdrop:a} $D$ contains an isometric copy of a disk
whose boundary has curvature $\kappa$;

\item\label{thickdrop:b} if the length of $\partial D$ is less than the length
of the circle with curvature $\kappa$
on the unit sphere then $D$ contains an isometric copy of a unit half-sphere.

\end{enumerate}

\end{thm}

\parit{Proof; (\ref{thickdrop:a}).} Recall that the {\it cut locus}
of $D$ with respect to its boundary $\partial D$ is defined as the closure
of the set of all points $x\in D$ such that the restriction of the distance function $\dist_x|_{\partial D}$
attains its global minimum at two or more points of $\partial D$.
The cut locus will be denoted as $\CutLoc D$.

After a small perturbation
of $\partial D$ we may assume that
$\CutLoc D$ is a graph embedded in
$D$ with finite number of edges.

\begin{wrapfigure}{r}{45mm}
\begin{lpic}[t(-0mm),b(0mm),r(0mm),l(0mm)]{pics/drop(1)}
\lbl[b]{22,48;$z$}
\lbl[lb]{31,52;$\bar z$}
\lbl[lt]{22,17;$y$}
\lbl[lb]{24,22;$\bar y$}
\end{lpic}
\end{wrapfigure}

Note that $\CutLoc D$ is a
deformation retract of $D$.
The retraction can be obtained by moving each point $y\in D\setminus \CutLoc D$ towards $\CutLoc D$
along the geodesic containing $y$ and
the point $\bar y\in\partial D$ closest to $y$.
In particular $\CutLoc D$ is a tree.

Since $\CutLoc D$  is a tree, it has
at least two vertices of valence one.
Among all points of $\partial D$ only
the non-smooth point of $\partial D$ with angle less
than $\pi$ belongs to $\CutLoc D$.
So there is at least one point $z$
of $\CutLoc D$ of valence one
contained in the interior of $D$.
The point $z$ has to be a focal point of $\partial D$;
this means that
the disk of radius $\dist_{\partial D}z$ centered at $z$ touches $\partial D$
with multiplicity at least two at some point $\bar z$.
At $\bar z$ the curvature of
the boundary of the disk centred at $z$ equals the curvature of
$\partial D$  and so it is at most $\kappa$. 
So this disk contains a disk with boundary of curvature $\kappa$.

\parit{(\ref{thickdrop:b}).}
By (\ref{thickdrop:a}) we can assume that $\kappa>0$.
Consider a locally isometric
immersion  of $D$ into the unit sphere, $\phi: D\looparrowright \mathbb{S}^2$.
Since the length of $\partial D$ is less than
$2{\cdot}\pi$, by Crofton's formula,  $\partial D$
does not intersect one of equators.
Therefore the curve $\phi(\partial D)$ is contained in a half-sphere, say $\mathbb{S}^2_+$.

Note that
it is sufficient to show that $\phi(D)$ contains the complement of $\mathbb{S}^2_+$.
Suppose the contrary; note that in this case $\phi(D)\subset \mathbb{S}^2_+$.
Applying (\ref{thickdrop:a}),
we get that $\phi(D)$ contains a disc bounded by a circle,
say $\sigma_\kappa$, of curvature $\kappa$.
Note that $\partial[\phi(D)]$ cuts $\sigma_\kappa$ from its antipodal circle;
therefore
$$\length \partial[\phi(D)]\ge \length \sigma_\kappa.$$
Note that
$$\length\partial D\ge \length \partial[\phi(D)].$$
On the other hand, by the assumptions
$$\length\partial D<\length \sigma_\kappa,$$
a contradiction.
\qeds

\begin{thm}{Lemma}\label{twocases}
Assume that $\Sigma$ is a spherical polyhedral $3$-manifold
with an isometric $\RR$-action, whose orbits are geodesic.
Then the quotient $\Lambda=(\Ram\Sigma)/\RR$ is a spherical polyhedral surface
of curvature $4$, and there are two possibilities.
\begin{enumerate}[(a)]
\item\label{twocases-a}
If $\Lambda$ is not contractible then it is isometric to the sphere of curvature $4$, further denoted as $\tfrac12\cdot\mathbb{S}^2$.
In this case, $\Ram \Sigma$ is isometric to the unit $\mathbb{S}^3$ or to $\Ram_{\mathbb{S}^1}\mathbb{S}^3$, where  $\mathbb{S}^1$ is a closed geodesic in $\mathbb{S}^3$.
\item\label{twocases-b} If $\Lambda$ is contractible
then a point $x\in \Ram \Sigma$ is singular
if and only if so is its projection $\bar x\in \Lambda$.
Moreover, the angle around each singular point $\bar x\in \Lambda$
is infinite.
\end{enumerate}
\end{thm}

\parit{Proof.}
We will consider two cases.

\parit{Case 1.} Assume the action $\RR\acts\Sigma$ is not periodic;
that is, it does not factor through an $\mathbb{S}^1$-action.
Then the group of isometries of $\Sigma$ contains a torus $\TT^2$.

From \cite[Proposition 3.9]{panov} one can deduce that the Euclidean
cone over $\Ram\Sigma$ is isometric to the ramification of $\RR^4$
in one $2$-plane or in a pair of two orthogonal $2$-planes.
It follows that $\Ram \Sigma$ is either
$$\Ram_{\mathbb{S}^1}\mathbb{S}^3 \ \ \mathrm{or}\ \ \Ram_{\mathbb{S}_a^1\cup \mathbb{S}_b^1}\mathbb{S}^3$$
where $\mathbb{S}_a^1$ and $\mathbb{S}_b^1$ are two opposite Hopf circles.
In both cases, the $\RR$-action is lifted from the Hopf $\mathbb{S}^1$-action on $\mathbb{S}^3$.

If $\Ram \Sigma=\Ram_{\mathbb{S}^1}\mathbb{S}^3$ then
$\Lambda=\tfrac12\cdot\mathbb{S}^2$ and therefore (\ref{twocases-a}) holds.

If $\Ram \Sigma=\Ram_{\mathbb{S}_a^1\cup \mathbb{S}_b^1}\mathbb{S}^3$ then  $\Lambda=\Ram_{\{a,b\}}(\tfrac12\cdot\mathbb{S}^2)$ where $a$ and $b$ are two poles of the sphere;
therefore (\ref{twocases-b}) holds.

\parit{Case 2.} Assume that the $\RR$-action is periodic.
Let $s$ be the number of orbits in the singular locus $\Sigma^{{\star}}$
and let $m$ be the number of multiple orbits in the regular locus $\Sigma^\circ$.

Note that the space $\Sigma^\circ/\mathbb{S}^1$ is an orbifold with constant curvature $4$;
it has $m$ orbi-points.
Passing to the completion of $\Sigma^\circ/\mathbb{S}^1$, we get $\Sigma/\mathbb{S}^1$.
In this way we add $s$ points to $\Sigma^\circ/\mathbb{S}^1$
which we will call the \emph{punctures};
this is a finite set of points formed by the projection of the singular locus $\Sigma^{{\star}}$
in the quotient space $\Sigma/\mathbb{S}^1$.

Now we will consider a few subcases.

Assume $s=0$; in other words $\Sigma^{{\star}}=\emptyset$.
Then $\Ram\Sigma$ is isometric to $\mathbb{S}^3$
and the $\RR$-action factors through the standard Hopf action;
that is, the first part of (\ref{twocases-a}) holds.

Assume either $s\ge 2$ or $s\ge 1$ and $m\ge 2$.
Then the orbifold fundamental group of $\Sigma^\circ/\mathbb{S}^1$ is infinite, the universal orbi-cover is a disk and it branches
infinitely over every puncture of $\Sigma/\mathbb{S}^1$.
The completion of the cover is contractible; 
that is, (\ref{twocases-b}) holds.

It remains to consider the subcase $s=1$ and $m=1$.
In this subcase the universal orbi-cover of $\Sigma^\circ/\mathbb{S}^1$ is a once-punctured $\mathbb{S}^2$ of curvature $4$ and $\Ram \Sigma=\Ram_{\mathbb{S}^1}\mathbb{S}^3$; 
that is, (\ref{twocases-a}) holds.
\qeds

\parit{Proof of Theorem \ref{thm:3D-sphere}.}
Suppose first $\Lambda=(\Ram \Sigma)/\RR$ is not contractible.
By Lemma \ref{twocases}, $\Lambda$ is isometric to $\mathbb{S}^2$,
and the ramification $\Ram\Sigma$ is
isometric either to $\mathbb{S}^3$ or $\Ram_{\mathbb{S}^1}\mathbb{S}^3$.
Both of these spaces are $\CAT[1]$;
so the theorem follows.

From now on we consider the case when $\Lambda$ is contractible and will
prove in this case that $\Ram\Sigma\in \CAT[1]$
if and only if $\Sigma$ is $\frac{\pi}{4}$-extendable.

\parit{If part.}
From Lemma \ref{twocases} it follows that $\Ram \Sigma$ branches infinitely
over singular circles of $\Sigma$.
So  $\Ram \Sigma$ is locally $\CAT[1]$ and we only need to show that
any closed geodesic $\gamma$ in $\Ram \Sigma$
has length at least $2{\cdot}\pi$
(see Proposition \ref{geodesic2pi}).

Let $\gamma$ be a closed geodesic in $\Ram \Sigma$;
denote by $\bar\gamma$ its projection in $\Lambda$.
The curve $\bar\gamma$ is composed
of arcs of constant curvature, say $\kappa$, joining  singularities of $\Lambda$.
Moreover for each singular point $p$ of $\Lambda$ that belongs to $\bar\gamma$ the angle
between the arcs of $\bar\gamma$ at $p$ is at least $\pi$.

Both of the above statements are easy to check;
the first one is also proved in \cite[Lemma 3.1]{panov1}.
The following lemma follows directly from \cite[Proposition 3.6 2)]{panov1}.

\begin{thm}{Lemma}\label{selfint} 
Assume $\Sigma$ and $\Lambda$ are as in the formulation of Lemma~\ref{twocases}.
Then for every geodesic $\gamma$
in $\Ram \Sigma$ its projection $\bar\gamma$ in $\Lambda$ has
a point of self-intersections.
\end{thm}

Summarizing all the above,
we can choose two sub loops in $\bar\gamma$,
say $\bar\gamma_1$ and $\bar\gamma_2$, which bound disks on $\Ram \Sigma/\mathbb R$
and
both of these disks satisfy the conditions of Lemma  \ref{thickdrop} for some $\kappa$.
Clearly, we can chose $\bar\gamma_1$ and $\bar\gamma_2$ so that $\bar\gamma_1\cap \bar\gamma_2$
is at most a finite set.

By our assumptions the disks bounded by $\bar\gamma_i$ cannot contain points on distance more than
$\frac{\pi}{4}$ from their boundary, otherwise $\Sigma$
would not be $\frac{\pi}{4}$-extendable.
So  we deduce from Lemma \ref{thickdrop}(\ref{thickdrop:b})
that
$$\length\bar\gamma_i\ge \ell(\kappa),\eqno({*})$$
where $\ell(\kappa)=\tfrac{2\cdot\pi}{\sqrt{\kappa^2+4}}$ is the length of a circle of curvature $\kappa$ on the sphere of radius $\tfrac12$.

Let $\alpha$ be an arc of $\gamma$
and $\bar\alpha$ be its projection in $\Lambda$.
Note that
$$\length \alpha
=
\tfrac{\pi}{\ell(\kappa)}\cdot\length \bar\alpha.
$$
Together with $({*})$,
this implies that $\length\gamma\ge 2{\cdot}\pi$.

\parit{Only if part.}
Suppose now that $\Sigma$ contains an immersed copy
of a ball with radius $\frac{\pi}{4}+\varepsilon$.
Consider a lift of this ball to $\Ram \Sigma$ and denote it by $B$.

Set as before $\Lambda=(\Ram \Sigma)/\RR$.
The projection of $B$ in $\Lambda$ is a
disc, say $D$,
of radius $\frac{\pi}{4}+\varepsilon$ and curvature $4$, isometrically
{\it immersed} in $\Lambda$. Since $\Lambda$ is contractible $D$
has to be {\it embedded} in $\Lambda$.

Consider a closed geodesic
$\bar \gamma\subset  \Lambda\setminus D$ which is obtained from $\partial D$ by a curve-shortening process.
Such a geodesic has to contain at least two singular points;
let $x$ be one of such points.
Choose now a lift  of $\bar\gamma$ to a horizontal geodesic path $\gamma$
on  $\Ram \Sigma$ with two (possibly distinct) ends at the $\RR$-orbit over $x$.

Finally consider a deck transformation
$\iota$ of $\Ram \Sigma$ that fixes  the $\RR$-orbit over $x$
and rotates around it $\Ram \Sigma$ by an angle larger than $\pi$.
The union of $\gamma$ with $\iota\circ\gamma$ forms a closed geodesic in
$\Ram \Sigma$ of length less than $2{\cdot}\pi$.
\qeds

The following statement is proved by the same methods as in the theorem.

\begin{thm}{Corollary from the proof}\label{linesinC2} Let $n\ge 2$ be an integer
and $X$ be a union of $n$ fibers
 of the Hopf fibration on  the unit $\mathbb{S}^3$.
Then $\Ram_{X}\mathbb{S}^3$ is $\CAT[1]$ if and only if there is no point on $\mathbb{S}^3$ at distance more than $\frac{\pi}{4}$ from $X$.
\end{thm}

\section{Case $(\CC^2)$} \label{pkspaces}

\parit{Proof of Theorem \ref{thm:main}; case $(\CC^2)$.}
Let us show that $\Ram_{\Gamma}\CC^2$ is $\CAT[0]$.

First assume that $\Gamma$ is finite.
Without loss of generality, we can assume that the origin is fixed by $\Gamma$.
Let $L$ be the union of all the lines in $\CC^2$
fixed by some non-identity elements of $\Gamma$.
Note that $\Ram_{\Gamma}=\Ram_{L}\CC^2$;
here $\Ram_AX$ denotes the completion of the universal cover of $X\backslash A$.

If $L=\emptyset$
or $L$ is a single line,
the statement is clear.

Set $\Theta=\mathbb{S}^3\cap L$; this
is a union of Hopf circles.
If the circles in $\Theta$
satisfy the conditions of
Corollary \ref{linesinC2} then $\Ram_{\Theta} \mathbb{S}^3$ is $\CAT[1]$.
Therefore
\[\Ram_\Gamma=\Ram_L\CC^2=\Cone(\Ram_{\Theta} \mathbb{S}^3)\in \CAT[0].\]

Suppose now that the conditions of
Corollary \ref{linesinC2} are not satisfied.
Denote by $\Xi$ the projection of $\Theta$ in $\tfrac12\cdot\mathbb{S}^2=\mathbb{S}^3/\mathbb{S}^1$;
note that $\Xi$ is a finite set of points.
In this case there is an open half-sphere containing all points $\Xi$.
Denote by $P$ the convex hull of $\Xi$.
Note that $\Xi$ and therefore $P$ are $\Gamma$-invariant sets.
Therefore the action on $\mathbb{S}^2$ is cyclic.
The latter means that $L$ consists of one line.

If $\Gamma$ is infinite, we can apply the above argument to each isotropy group
of $\Gamma$.
We get that $\Ram_{\Gamma_x}$ is $\CAT[0]$
for the isotropy group $\Gamma_x$ at any point $x\in\CC^2$.
Then it remains to apply Proposition \ref{orbianalog}(\ref{orbianalog:i}).

Now let us show that the inclusion $W_\Gamma\hookrightarrow \Ram_\Gamma$ is a homotopy equivalence.
Fix a singular point $y$ in $\Ram_\Gamma$
and let $x$ be its projection to $\CC^2/\Gamma$.
By Proposition \ref{orbianalog}(\ref{orbianalog:ii})
the link at $y$ is the same as the link of the ramification of the cone at $x$.
The latter space is the ramification of $\mathbb{S}^3$ in a non-empty collection of Hopf circles,
which is clearly contractible.
It remains to apply Allcock's lemma \ref{lem:allcock}.
\qeds

\section{The counterexample}

In this section we use the technique introduced above to show that
the answer to the Question \ref{quest:generalized} is negative
without additional assumptions on $\mathcal{P}$.

\begin{thm}{Theorem}\label{notCAT[1]} There is a positively curved spherical polyhedral space $\mathcal{P}$ 
such that $\Ram \mathcal{P}$ is not $\CAT[1]$.
Moreover, one can assume that $\mathcal{P}$ is homeomorphic to $\mathbb{S}^3$
and it admits an isometric $\mathbb{S}^1$-action
with geodesic orbits.
\end{thm}

\parit{Proof.} 
Consider a triangle $\Delta$ on the sphere of curvature $4$ with one angle $\frac{\pi}{n}$ and the other two $\frac{\pi{\cdot}(n+1)}{2{\cdot} n}\z+\varepsilon$;
here $n$ is a positive integer and $\eps>0$.
Note that two sides of $\Delta$
are longer than $\frac{\pi}{4}$.

Denote by $\Lambda$ the doubling of $\Delta$.
The space $\Lambda$ is a spherical polyhedral space with curvature $4$;
it has three singular points which correspond to the vertices of $\Delta$.
Label the point with angle $\frac{2\cdot\pi}{n}$ by $x$.

According to \cite[Theorem 1.8]{panov} there is a unique up to isometry
polyhedral spherical space $\mathcal{P}$
with an isometric action $\mathbb{S}^1\acts\mathcal{P}$
such that $\mathbb{S}^1$-orbits are geodesic,
$\Lambda$ is isometric to the quotient space $\mathcal{P}/\mathbb{S}^1$
and the point $x$ corresponds to the orbit of multiplicity $n$,
while the rest of orbits are simple.

Note that the points in $\mathcal{P}$ 
on the $\mathbb{S}^1$-fiber over $x$ are regular.
The distance from this fiber to the singularities of
$\mathbb{S}^3$ is more than $\frac{\pi}{4}$;
that is, $\mathcal{P}$ is not $\frac{\pi}{4}$-extendable.
By Theorem \ref{thm:3D-sphere} we conclude that $\Ram \mathcal{P}$ is not $\CAT[1]$.
\qeds

\section{Line arrangements}\label{sec:linearrangement}

The following theorem is the main result of this section.

\begin{thm}{Theorem}\label{linearrangement} Let $\mathcal{P}$ be a non-negatively curved polyhedral space homeomorphic to $\CP^2$ whose singularities form a complex line arrangement on $\CP^2$.
Then $\Ram \mathcal{P}$ is a $\CAT[0]$ space and the inclusion $(\Ram \mathcal{P})^\circ\hookrightarrow\Ram \mathcal{P}$ is a homotopy equivalence.
\end{thm}

It follows that all complex line arrangements in $\CP^2$
appearing as singularities of non-negatively curved polyhedral metrics
have aspherical complements. The class of such arrangements is characterized in
Theorem \ref{general};
this class includes all the arrangements from Theorem \ref{cor:3n}.

\parit{Proof.}
According to \cite{cheeger}
and \cite{panov}
the metric on $\mathcal{P}$ is
polyhedral K\"ahler.

First let us show that $\Ram \mathcal P$ is $\CAT[0]$.
By Theorem \ref{orbianalog},
it is sufficient to show that the ramification of the cone of each singular point $x$ in $\mathcal P$ is $\CAT[0]$.

If there are exactly two lines meeting at $x$ then the cone of $x$ is a direct product by Theorem \ref{allPK}(\ref{allPK:3}),
and the statement is clear.

If  more than two lines meet at $x$ consider the link $\Sigma$
of the cone at $x$. According to Theorem \ref{allPK} there is a free $\mathbb{S}^1$-action on $\Sigma$
inducing on it the structure of the Hopf fibration.
The quotient $\Sigma/\mathbb{S}^1$ is a $2$-sphere
with spherical polyhedral metric of curvature $4$
and the conical angle is at most $2\cdot\pi$
around any point.
It follows from Zalgaller's theorem
that $\Sigma$ is $\frac{\pi}{4}$-extendable. So by Theorem \ref{thm:3D-sphere}
$\Ram \Sigma$ is $\CAT[1]$.

It remains to show that $\Ram \mathcal P^{\circ}\hookrightarrow \Ram \mathcal P$ is a homotopy equivalence.
The latter follows from Allcock's lemma \ref{lem:allcock}
the same way as at the end of the proof of Theorem  \ref{thm:main}, case $\CC^2$.
\qeds

\parit{Proof of Theorem \ref{cor:3n}.} By Theorem \ref {linearrangement}
it suffices to know that there is a
non-negatively curved polyhedral metric on $\CP^2$ with singularities
at the line arrangement. 
It is shown in \cite{panov} that for any arrangement of $3\cdot n$ lines that satisfies Hirzebruch's property and such that no $2\cdot n$ lines of the arrangement pass through one point,
such a metric exists. 

We are left with the case when at least $2\cdot n$ lines of the arrangement pass through one point, say $p$. Take any other line that does not pass through $p$. This line has at least $2\cdot n$ distinct intersections with other lines of the arrangement. So $n+1\ge 2\cdot n$, and we conclude that the arrangement is composed of three generic lines; hence it complement is aspherical.
\qeds

\parbf{General line arrangements.}
Let $(\ell_1,\dots,\ell_n)$ be a line arrangement in $\CP^2$.
The number of lines $\ell_i$ passing through a given point $x\in\CP^2$
will be called the \emph{multiplicity} of $x$, briefly
\def\mult{\mathrm{mult}}
$\mult_x$.

Let us associate to the arrangement a symmetric $n\times n$ matrix $(b_{ij})$.
For $i\ne j$
put $b_{ij}=1$ if the point $x_{ij}= \ell_i\cap \ell_j$
has multiplicity $2$
and $b_{ij}=0$ if its multiplicity is $3$ or higher.
The number $b_{jj}+1$ equals the number of points on $\ell_j$ with the multiplicity $3$ and higher.

The next theorem follows from \cite[Theorem 1.12 and Lemma 7.9]{panov};
it reduces the existence of a non-negatively curved polyhedral K\"ahler metric on $\CP^2$
with singularities at a given
line arrangement to the existence of a solution of certain system of linear equalities and inequalities.

\begin{thm}{Theorem}\label{general} Let $(\ell_1,\dots,\ell_n)$ be a line arrangement in $\CP^2$ and $(b_{ij})$ be its matrix.
There exists a non-negatively curved polyhedral K\"ahler metric on $\CP^2$
with the singular locus formed by the lines $\ell_i$
if and only there are real numbers $(z_1,\dots,z_n)$
such that
\begin{enumerate}[(i)]
\item \label{general:i}for each $k$, we have
\[0<z_k<1;\]
\item \label{general:ii}for each $j$, we have
\[\sum_{k} b_{jk}\cdot z_k=1\]
and \[\sum_{k} z_k=3;\]
\item\label{general:iv} for each point $x\in\CP^2$ with multiplicity at least 3, we have%
\[\alpha_x=1-\tfrac{1}{2}\cdot\sum_{\{k|x\in \ell_k\}}z_k>0.\]
\end{enumerate}
\end{thm}

Let us explain the geometric meaning of the above conditions.
If $(z_1,\dots,z_n)$ satisfy the conditions then
there is a polyhedral K\"ahler metric on $\CP^2$
with the conical angle around $\ell_i$ equal to $2\cdot\pi\cdot(1-z_i)$.
The inequalities (\ref{general:i}) say that conical angles are positive and less than $2\pi$.

Each of $n$ equalities (\ref{general:ii}) 
is the Gauss--Bonnet formula
for the flat metric with conical singularities at a line of the arrangement;
the additional equality expresses the fact that the canonical bundle of
$\CP^2$ is $O(-3)$;
that is, the cube of the tautological line bundle.

The link $\Sigma_x$ at $x$ with the described metric
is isometric to a $3$-sphere with an $\mathbb{S}^1$-invariant metric.
A straightforward calculation shows that
the length of an $\mathbb{S}^1$-fiber in $\Sigma_x$ is $2\cdot\pi\cdot\alpha_x$,
where $\alpha_x$ as in (\ref{general:iv}).
Equivalently, $\pi\cdot\alpha_x$ is the area of the quotient space $\Sigma_x/\mathbb{S}^1$.

The construction of the metric in this theorem relies on a parabolic version of Kobayshi--Hitchin correspondence
established by Mochizuki \cite{mochizuki}. 
Surprisingly, the system of $n$ linear equations in
(\ref{general:ii}) is equivalent to the following quadratic equation. 
(The equation 
implies the system by  \cite[Lemma 7.9]{panov} and  the converse implication 
is a direct computation.)
$$\sum_{\{x|\mult_x>2\}}(\alpha_x-1)^2
-
\sum_{j=1}^nz_j^2\cdot b_{jj}
=
\tfrac{3}{2}.$$
This equation is the border case of a parabolic Bogomolov--Miayoka inequality. Geometrically
it expresses the second Chern class of $\CP^2$ as a sum of contributions of
singularities of the metric.

The following corollary generalizes Theorem \ref{cor:3n}.

\begin{thm}{Corollary}\label{generalarrangement}
Any line arrangement $(\ell_1,\dots,\ell_n)$ in $\CP^2$
for which one can find positive $z_j$ satisfying equalities and inequalities
of Theorem  \ref{general} has an aspherical complement.
\end{thm}

The arrangements of lines
as in Theorem \ref{cor:3n}
satisfy the conditions in Theorem~\ref{general}
with $z_i =\tfrac1n$
at all $3\cdot n$ lines of the arrangement. 
This is proved by an algebraic computation; 
see \cite[Corollary 7.8]{panov}.
The restriction that at most $2\cdot n-1$ lines pass through one point
follows from (\ref{general:iv}).
Therefore the corollary above is a generalization of  Theorem~\ref{general}.

\parit{Proof.}
By Theorem \ref{general} there is a non-negatively curved polyhedral metric on
$\CP^2$ with singularities along $(\ell_1,\dots,\ell_n)$ and so one can apply Theorem \ref{linearrangement}.
\qeds

\end{document}